\documentclass[letterpaper, 10pt, journal, twoside]{IEEEtran}
\IEEEoverridecommandlockouts
\usepackage{graphicx,keyval,trig}
\usepackage{amsfonts,amssymb,dsfont,epstopdf, epsfig}
\usepackage{amsthm}

\usepackage{enumitem}
\usepackage{float}
\usepackage{mathtools, mathrsfs, subfigure,amsmath}
\usepackage[usenames, dvipsnames]{color}
\usepackage[colorlinks=true,
         raiselinks=true,
         linkcolor=MidnightBlue,
         citecolor=Mahogany,
         urlcolor=ForestGreen,
         pdfauthor={},
         pdftitle={},
         pdfkeywords={stochastic reachability and avoidance, optimal
control},
         pdfsubject={Technical Report},
         plainpages=false]{hyperref}

  \usepackage{todonotes}

\usepackage{makecell}

\allowdisplaybreaks

\usepackage{tikz}
\usetikzlibrary{chains}
\usetikzlibrary{fit}
\usepackage{pgflibraryarrows}		%optional
\usepackage{pgflibrarysnakes}		%optional
\usepackage{xcolor}
\usepackage{epsfig}
\usetikzlibrary{shapes.symbols,patterns} % for source symbols
\usepackage{pgfplots}

        \allowdisplaybreaks

\def \triple|{|\! | \! |}

\def\<{\langle}
\def\>{\rangle}
\def\~{\tilde}

\newcommand{\drv}{\ensuremath{\mathrm{d}}}

\usepackage{verbatim}
\newcommand{\cd}{C^1(X)}
\newcommand{\cds}{C^1(X)^*}
\newcommand{\con}{C(X\times U)}
\newcommand{\me}{\mathcal{M}(X\times U)}
\newcommand{\inner}[3]{\left\langle #1,#2\right\rangle_{#3}}

\newcommand{\norma}[2]{\left\|\,#1\,\right\|_{#2}}

\newcommand{\oloklirosi}[4]{\int_{#1}^{#2}#3\,\,\drv#4}

\newcommand{\athroisma}[4]{\sum_{#1=#2}^{#3}#4}

\newcommand{\ar}{\mathbb{R}}

\newtheorem{thm}{Theorem}[section]
\newtheorem{protasi}{Proposition}[section]
\newtheorem{orismos}{Definition}[section]
\newtheorem{porisma}{Corollary}[section]

\newtheorem{limma}{Lemma}[section]

\newtheorem{ass}{Assumption}[section]

\title{On Infinite Linear Programming and the Moment Approach to Deterministic Infinite Horizon Discounted Optimal Control Problems
\thanks{{Research was supported by the European Union 7th Framework Program "Scalable Proactive Event-Driven Decision-making (SPEEDD)".}}}
\author{
Angeliki Kamoutsi, Tobias Sutter, Peyman Mohajerin Esfahani, and John Lygeros
\thanks{The authors are with the Automatic Control
Laboratory, ETH Z\"urich, Switzerland and the Delft Center for Systems and Control, TU Delft, Netherlands; Emails: {\{kamoutsa, sutter, lygeros\}@control.ee.ethz.ch, P.MohajerinEsfahani@tudelft.nl}
}}

\begin{document}
     \maketitle
     \thispagestyle{empty}
\begin{abstract}
We revisit the linear programming approach to deterministic, continuous time, infinite horizon discounted optimal control problems. In the first part, we relax the original problem to an infinite-dimensional linear program over a measure space and prove equivalence of the two formulations under mild assumptions, significantly weaker than those found in the literature until now. The proof is based on duality theory and mollification techniques for constructing approximate smooth subsolutions to the associated Hamilton-Jacobi-Bellman equation. In the second part, we assume polynomial data and use Lasserre's hierarchy of primal-dual moment-sum-of-squares semidefinite relaxations to approximate the value function and design an approximate optimal feedback controller. We conclude with an illustrative example.
\end{abstract}
% Keywords appear just beneath the abstract. Use only for final version
\section{Introduction} \label{sec:introduction}
% Drop letter for first word of the Introduction
% Here we have the typical use of a "T" for an initial drop letter
% and "HIS" in caps to complete the first word.
% Use only for final version
\IEEEPARstart{W}{e} study the infinite horizon optimal control problem (OCP) with discounted payoff subject to state constraints.
A comprehensive theoretical framework has been developed over the years to tackle such OCPs, e.g., via the Pontryagin's maximum principle and the associated Hamilton-Jacobi-Bellman equation (HJB)~\cite[Sections III.2, IV.5]{Bardi}. However, these sophisticated and mathematically elegant methods are not always easy to implement in practice. Several approaches, including shooting methods~\cite{shooting}, model predictive control~\cite{MPC}, direct methods~\cite{direct} and neural-network-based algorithms~\cite{NN}, have been proposed in the literature to alleviate this difficulty.

In this work, we study an approximate dynamic programming approach based on infinite-dimensional linear programming (LP). The LP formulation of deterministic and stochastic optimal control problems has been studied extensively in the literature~\cite{HLL2, Gaitsgory, Vinter, BhattBorkar, FlemingVermes}. The main idea is to embed the set of admissible controls for the original problem in a new space, a space of measures. In our particular setting, we identify each admissible control with a discounted occupation measure and observe that it satisfies a linear equation that provides a measure-theoretic interpretation of the system dynamics. In this way, we introduce an infinite-dimensional LP (primal) over an appropriate space of measures. The optimal value of the primal LP is no greater than the optimal value of the OCP and, under mild assumptions described in Theorem~\ref{equivalence}, the two are equal. Using duality theory, we introduce the dual LP, which involves finding the supremum of all smooth subsolutions to the associated Hamilton-Jacobi-Bellman (HJB) equation, and prove that there is no duality gap. As a result, we derive a characterization of the value function as the upper envelope of the smooth subsolutions to the HJB equation. It is worth mentioning that the LP approach is particularly appealing for dealing with unconventional problems involving additional constraints or secondary costs, where traditional dynamic programming techniques are not applicable. In particular, in the OCP of our interest, the state and input constraints are directly incorporated into the measure space associated with the primal program and the function space associated with the dual program.	

We then, use finite-dimensional approximations of the primal LP and its dual to approximate the optimal value function and extract a near optimal feedback control. In the literature, several works propose finite LP approximations based on state-and-control-space gridding \cite{Gaitsgory,ref:Peyman-17}. On the other hand, more recently, Lasserre et al. in~\cite{LasserreHenrion, controlsynthesis}, used convex optimization methods and proposed a hierarchy of finite-dimensional semidefinite programming (SDP) relaxations to tackle nonlinear finite horizon OCPs with polynomial data. We extend the latter approach to the infinite horizon case.
%===========================================================================================================================================================

\textbf{Contribution.}
This paper is inspired and extends the techniques discussed in~\cite{LasserreHenrion},~\cite{controlsynthesis} and ~\cite{Vinter}. We show equivalence of the LP formulation to the OCP (Theorem~\ref{equivalence}) under weaker assumptions than those known in the literature~\cite{Gaitsgory}; in particular, our assumptions do not require Lipschitz continuity of $V_C^{Y_{\delta}}$ with respect to $\delta$~\cite[Theorem 4.4]{Gaitsgory}, a technical assumption that usually cannot be checked even for simple OCPs. Moreover Assumption 1 in~\cite{Gaitsgory}, needs to hold only for the cost function of the OCP and not for all continuous functions. Our proof of equivalence is based on duality theory and mollification techniques for constructing approximate smooth subsolutions to the HJB equation, also proposed in~\cite{Vinter} for finite horizon differential inclusions and in~\cite{FlemingVermes} for controlled diffusions. Our finite-dimensional optimization approximations are based on Lasserre's hierarchy of semidefinite relaxations~\cite{LasserreHenrion} instead of the grid-based LP approximations proposed in~\cite{Gaitsgory}. However, unlike~\cite{LasserreHenrion}, we treat the infinite horizon polynomial OCP by introducing the notion of discounted occupation measures. Finally, recently Korda et al.~\cite{Korda} have presented a new approach, where the primal LP is tightened and optimization is carried out over polynomial densities of measures. In this way, they are able to provide a controller design with strong convergence guarantees. However, their initial assumptions restrict significantly the class of OCPs for which their method is applicable. Although in our controller extraction, which follows closely the reasoning in~\cite{controlsynthesis}, we do not provide convergence guarantees, our simulation results show better performance in the presented numerical example.
%===============================================================================
%===============================================================================
%===============================================================================

\textbf{Basic definitions and notations.}
Let $\mathbb{X}$ be a normed vector space. The \emph{weak* topology} on its topological dual $\mathbb{X}^*$ is the smallest topology with respect to which the linear functionals in $\{\mathbb{X}^*\ni x^*\rightarrow \inner{x^*}{x}{}:=x^*(x){}: x\in \mathbb{X}\}$ are continuous. Let $\mathcal{X}$ be a metrizable  topological space and let $\mathcal{B}(\mathcal{X})$ be its Borel $\sigma$-algebra. Let $C(\mathcal{X})$ be the Banach space of real-valued bounded continuous functions on $\mathcal{X}$ together with the sup-norm $\norma{}{\infty}$. We denote by $\mathcal{M}(\mathcal{X})$ the Banach space of finite signed Borel measures on $\mathcal{X}$ equipped with the total variation norm and by $\mathcal{P}(\mathcal{X})$ the set of Borel probability measures. Let $\delta_x\in\mathcal{P}(\mathcal{X})$ be the Dirac measure centred on $x\in\mathcal{X}$. When $\mathcal{X}$ is compact,  $C(\mathcal{X})^*\simeq\mathcal{M}(\mathcal{X})$. In particular, we identify each $\mu\in\mathcal{M}(\mathcal{X})$ with the bounded linear functional $\inner{\mu}{\cdot}{}:C(\mathcal{X})\rightarrow\ar$,
$
  \inner{\mu}{l}{}:=\oloklirosi{\mathcal{X}}{}{l}{\mu},
$
for all $l\in C(\mathcal{X})$. Moreover, if $\mathcal{M}(\mathcal{X})^+$ is the convex cone of finite nonnegative Borel measures on $\mathcal{X}$, then its dual convex cone is the set $C(\mathcal{X})^+$ of nonnegative continuous functions on $\mathcal{X}$. The \emph{support} of a measure $\mu\in\mathcal{M}(\mathcal{X})^+$ is the smallest closed set whose complement has zero measure and is denoted by spt $\mu$. If $X\subset\ar^n$, we consider the Banach space
$
C^1(X)=\{\phi\in C(X):\,\,\,\frac{\partial\phi}{\partial x_i}\in C(X),\,\,\,\mbox{for all}\,\,\,i=1,\ldots,n\},
$
with
$
\norma{\phi}{\infty}^1=\norma{\phi}{\infty}+\athroisma{i}{1}{n}{\norma{\frac{\partial\phi}{\partial x_i}}{\infty}}.
$
Let $\ar[x]$ be the space of polynomials of the variable $x\in\ar^n$. Let $\alpha=(\alpha_1,\ldots,\alpha_n)^\top\in\mathbb{N}^n$ be a multi-index. A monomial is defined as
$
x^{\alpha}:=x_1^{\alpha_1}x_2^{\alpha_2}\ldots x_n^{\alpha_n}
$
and its degree is
$
|\alpha|:=\athroisma{i}{1}{n}{\alpha_i}.
$
 A polynomial $p\in\ar[x]$ of degree $\deg(p)=d$ is given by
$
p(x):=\sum_{|\alpha|\le d}p_{\alpha}x^{\alpha}.
$
Similarly, let $\ar[x,u]$ be the space of polynomials of the variable $(x,u)\in\ar^n\times\ar^m$. Then a polynomial $q\in\ar[x,u]$ of degree $d$ has the form
$q(x,u):=\sum_{\substack{\gamma\in\mathbb{N}^n\times\mathbb{N}^m\\|\gamma|\le d}}q_\gamma\,\,\, (x,u)^{\gamma}=\sum_{\substack{\alpha\in\mathbb{N}^n,\,\beta\in\mathbb{N}^m\\|\alpha|+|\beta|\le d}}q_{\alpha\beta}\,\,x^{\alpha}u^{\beta}.
$ Let $X\subset\ar^n$ and $Y\subset\ar^m$ be compact subsets. Given a measure $\mu\in\me^+$, the real number $z_{\gamma}:=\oloklirosi{X\times U}{}{(x,u)^{\gamma}}{\mu(x,u)}=\oloklirosi{X\times U}{}{x^{\alpha}u^{\beta}}{\mu(x,u)}$, is called its \emph{moment} of order $\gamma=(\alpha,\beta)\in\mathbb{N}^n\times\mathbb{N}^m$. Given a real sequence $z=(z_\gamma)_{\gamma\in\mathbb{N}^n\times\mathbb{N}^m}$, the \emph{Riesz linear functional} $L_z:\ar[x,u]\rightarrow\ar$ is defined by $L_z(q):=\sum_{\gamma\in\mathbb{N}^n\times\mathbb{N}^m} q_{\gamma}z_{\gamma}$, for each polynomial $q(x,u)=\sum_{\gamma\in\mathbb{N}^n\times\mathbb{N}^m} q_{\gamma}\,\,(x,u)^{\gamma}$. %Note that when the sequence $z$ has a representing measure $\mu\in\me^+$, then $\inner{\mu}{q}{}=L_z(q)$.
 Let $\nu_d(x,u):=((x,u)^{\gamma})_{|\gamma|\le d}$. The \emph{moment matrix} of order $d$ is given by $M_d(z):=L_z(\nu_d(x,u)\,\nu_d(x,u)^\top)$, where the latter notation means that $L_z$ is applied entrywise on the input matrix. Similarly, given $h\in\ar[x,u]$ the \emph{localizing matrix} of order $d$ w.r.t. $z$ and $h$ is defined by $M_d(h\,z):=L_z(h(x,u)\,\nu_d(x,u)\,\nu_d(x,u)^\top)$. A polynomial $h\in\ar[x,u]$ is a \emph{sum of squares (s.o.s.)} if it can be written as $h(x,u)=\sum_{j\in J}h_j^2(x,u)$, for some finite index set $J$, and $\{h_j:\,j\in J\}\subset\ar[x,u]$. We denote by $\Sigma[x,u]\subset\ar[x,u]$, the space of s.o.s. polynomials.
%===============================================================================
\section{Problem Statement} \label{sec:ProbStatement}
%===============================================================================
Let $X\subset\ar^n$ and $U\subset\ar^m$. We consider the following infinite horizon control problem:
\begin{equation}\label{B.1}
\begin{aligned}
% \nonumber to remove numbering (before each equation)
  \dot{x}(t) &= f(x(t),u(t)),\,\,\,\,\,\,\,t>0 \\
  x(0) &= x_0,
  \end{aligned}
\end{equation}
where $x_0\in X$ and $f:\ar^n\times U\rightarrow\ar^n$. A \emph{control} is a Borel measurable function $u(\cdot):[0,\infty)\rightarrow U$. The set of controls is denoted by $\mathcal{U}^0$. 	Under Assumption~\ref{ass:control:model}~\ref{ass:H1} and~\ref{ass:H2} below, for each $u(\cdot)\in\mathcal{U}^0$, the control system~(\ref{B.1}) admits a unique solution $x(\cdot|x_0,u(\cdot))$.  We impose the state constraint
\begin{equation}\label{statecon}
  x(t|x_0,u(\cdot))\in X,\,\,\,\,\,\,\,\mbox{for all}\,\,\,\,t\geq 0.
\end{equation}
The set of \emph{admissible controls} is given by
$
  \mathcal{U}_X(x_0):=\{u(\cdot)\in\mathcal{U}^0: \mbox{~(\ref{statecon}) holds}\}.
$
The OCP is described by
\begin{equation}\label{B.7}
V_X(x_0):=\inf_{u(\cdot)\in\mathcal{U}_X(x_0)} \oloklirosi{0}{\infty}{e^{-\lambda t}g(x(t|x_0,u(\cdot)),u(t))}{t},
\end{equation}
where $\lambda>0$ is the \emph{discount factor} and $g:\ar^n\times U\rightarrow\ar$. Consider the following conditions:
\begin{ass}[Control model]\label{ass:control:model}   \
\begin{enumerate}[label=\text{(H\arabic*)}, itemsep = -1mm, topsep = -3mm]
  \item \label{ass:H1} $f,g$ are bounded and continuous, Lipschitz in $x$ uniformly in $u$;
  \item \label{ass:H2} $X$ and $U$ are compact;
  \item \label{ass:H3} $\lambda>L$, where $L$ is the Lipschitz constant for $f$;
  \item $\mathcal{U}_X(x_0)\ne\emptyset$.
\end{enumerate}
\end{ass}
Under suitable assumptions \cite[Thm.~5.10]{Bardi}, the value function $V_X$ is the unique constrained viscosity solution of the \emph{Hamilton-Jacobi-Bellman (HJB) equation}
\begin{equation*}
  \lambda \phi(x)+\sup_{u\in U}\{-f(x,u)\cdot \nabla_x\phi(x) -g(x,u)\}=0,\,\,x\in\textup{Int}X.
\end{equation*}
In general, the infimum in~(\ref{B.7}) is not attained, so our next step is to consider a relaxation of the OCP. A \emph{relaxed control} or a \emph{Young measure} is a measurable map $m(\cdot):[0,\infty)\rightarrow\mathcal{P}(U)$, i.e., for each $h\in C(U)$, the scalar-valued function $[0,+\infty)\ni t\rightarrow\oloklirosi{U}{}{h(u)}{m_t(du)}$ is Lebesgue measurable. The set of relaxed controls is denoted by $\mathcal{V}^0$.  We consider the relaxed control system
\begin{equation}\label{relaxedsystem}
\begin{aligned}
% \nonumber to remove numbering (before each equation)
   \dot{x}(t)&=\oloklirosi{U}{}{f(x(t),u)}{m_t(u)},\,\,\,t\geq 0, \\
  x(0) &= x_0.
  \end{aligned}
\end{equation}
Given a relaxed control $m(\cdot)$, we denote by $x(\cdot|x_0,m(\cdot))$ the solution of the system~(\ref{relaxedsystem}). We impose the state constraint
\begin{equation}\label{relaxedconstraint}
  x(t|x_0,m(\cdot))\in X,\,\,\,\mbox{for all}\,\, t\geq 0.
\end{equation}
The set of admissible relaxed controls is given by $\mathcal{V}_X(x_0):=\{m(\cdot)\in\mathcal{V}^0: \mbox{~(\ref{relaxedconstraint}) holds}\}$. Note that $\mathcal{U}_X(x_0)\subset\mathcal{V}_X(x_0)$ by identifying each control $u(\cdot)$ with the Dirac measure-valued function $t\rightarrow \delta_{u(t)}$. The relaxed OCP is described by
\begin{equation}\label{relaxed OCP}
  V_X^R(x_0)\!:=\!\!\!\!\min_{m(\cdot)\in\mathcal{V}(x_0)}\oloklirosi{0}{\infty}{\!\!\!e^{-\lambda t}\oloklirosi{U}{}{\!g(x(t),u)}{m_t(u)}}{t}.
\end{equation}
By using the weak* compactness of $\mathcal{P}(X)$ one can prove that the minimum is indeed attained. Moreover, $V_X^R(x_0)\le V_X(x_0)$. We add the following assumption.
\begin{ass}[No relaxation gap]\label{norelaxationgap}
$V_X(x_0)=V_X^R(x_0)$, for all $x_0\in X$ such that $\mathcal{U}_X(x_0)\ne\emptyset$.
\end{ass}
Sufficient conditions under which Assumption~\ref{norelaxationgap} is satisfied are given e.g., by the assumptions of Filipov-Wazewski type theorems with state constraints~\cite{FW} or by conditions described in~\cite[Prop. 4.3]{Gaitsgory2}. In particular, Assumption~\ref{norelaxationgap} is satisfied in case of invariance of $X$ w.r.t. the solutions of the OCP or in case of uncontrolled dynamics. Morover in the case of input-affine dynamics, Assumption~\ref{norelaxationgap} is satisfied under Assumption~\ref{ass:control:model} and if $U$ is convex. For a related discussion on sufficient conditions see also the paragraph following Assumption 1 in~\cite{Gaitsgory2}.
%===============================================================================
%=======================================================================================================================================
\section{Linear programming Formulation} \label{sec:Method}
%======================================================================================================================================
For each $u(\cdot)\in\mathcal{U}_X(x_0)$, we define the discounted occupation measure $M^u\in\me^+$,
\begin{equation}\label{oom}
  M^u(E):=\oloklirosi{0}{\infty}{e^{-\lambda t}\delta_{(x(t|x_0,u(\cdot)),u(t))}(E)}{t},
\end{equation}
for all $E\in\mathcal{B}(X\times U)$. Moreover, for each $m(\cdot)\in\mathcal{V}_X(x_0)$, we define the discounted relaxed occupation measure $N^m\in\me^+$,
\begin{equation}\label{rom}
  N^m(E):=\!\oloklirosi{0}{\infty}{\!\!\!e^{-\lambda t}\oloklirosi{U}{}{\delta_{(x(t|x_0,m(\cdot)),u)}(E)}{m_t(u)}}{t},
\end{equation}
for all $E\in\mathcal{B}(X\times U)$.
Note that for every $l\in\con$,
\begin{align}
\inner{M^u}{l}{} &=\!\oloklirosi{0}{\infty}{e^{-\lambda t} \ l(x(t|x_0,u(\cdot)),u(t))}{t} \label{poom} \\
\inner{N^m}{l}{} &=\!\oloklirosi{0}{\infty}{\!\!\!\!e^{-\lambda t\!\!}\oloklirosi{U}{}{\!l(x(t|x_0,m(\cdot)),u)\!}{m_t(u)}}{t}. \label{prom}
\end{align}
By setting
$
\mathcal{M}_X:=\{M^u:u(\cdot)\in\mathcal{U}_X(x_0)\}
$
 and
 $
 \mathcal{R}_X:=\{N^m:m(\cdot)\in\mathcal{V}_X(x_0)\},
 $
we can rewrite the OCP~(\ref{B.7}) as
$
V_X(x_0)=\inf\{\inner{\mu}{g}{}:\mu\in\mathcal{M}_X\},
$
and the relaxed OCP~(\ref{relaxed OCP}) as
$
V_X^R(x_0)=\inf\{\inner{\mu}{g}{}:\mu\in\mathcal{R}_X\}.
$

In both cases, the cost is linear over $\mu$ but $\mathcal{R}_X$ and $\mathcal{M}_X$ are nonconvex sets. To transform the OCP into a convex program, we define the linear operator $A:C^1(X)\rightarrow C(X\times U)$ by
$
  (A\phi)(x,u):=\lambda\phi(x)-\nabla_x\phi(x)\cdot f(x,u),
$
for all $(x,u)\in X\times U,\,\phi\in C^1(X)$. Since $\norma{A\phi}{\infty}\le\left(\lambda+\sup_{(x,u)\in X\times U}|f(x,u)|\right)\norma{\phi}{\infty}^1,
$ for all $\phi\in C^1(X)$, $A$ is bounded and in particular it is weakly continuous. Let $A^*:\me\rightarrow\cds$ be the adjoint map of $A$, defined by $
\inner{A^*\mu}{\phi}{}=\inner{\mu}{A\phi}{}=\oloklirosi{X\times U}{}{A\phi}{\mu},
$ for all $\mu\in\me$ and $\phi\in\cd$. By the weak continuity of $A$, $A^*$ is well defined and weakly* continuous.
\begin{comment}
 In fact, one can check directly that $|\inner{A^*\mu}{\phi}{}|\le \left(\lambda+\sup_{(x,u)\in X\times U}|f(x,u)|\right)\norma{\phi}{\infty}^1\norma{\mu}{},$ for all $\mu\in\me$ and $\phi\in\cd$. So $A^*$ is a bounded linear operator with $\norma{A^*}{\mbox{op}}\le\lambda+\sup_{(x,u)\in X\times U}|f(x,u)|.$
\end{comment}
We introduce the linear program:
\begin{equation}\label{primalp}
  \mathcal{J}(x_0):=\inf\{\inner{\mu}{g}{}:\mu\in\mathcal{F}_X\},
\end{equation}
where,
$
{\mathcal{F}}_X:=\{\mu\in\me^+: A^*\mu=\delta_{x_0}\}
$
and $\delta_{x_0}$ is identified as an element of $\cds$. The following proposition establishes the relationship between the above programs and asserts that the infinite LP~(\ref{primalp}) admits a minimizer.
\begin{protasi}\label{minimizer}
  $\mathcal{M}_X\subset\mathcal{R}_X\subset\mathcal{F}_X$ and so $\mathcal{J}(x_0)\le V_X^R(x_0)\le V_X(x_0)$. Moreover, if $\mathcal{F}_X\ne\emptyset$, then the infimum in~(\ref{primalp}) is attained.
\end{protasi}
\begin{proof} Since $\mathcal{U}_X(x_0)\subset\mathcal{V}_X(x_0)$, we get $\mathcal{M}_X\subset\mathcal{R}_X$ and $V_X^R(x_0)\le V_X(x_0)$. Moreover, if $m(\cdot)\in\mathcal{V}(x_0)$, then for each $\phi\in\cd$ it holds that $\nabla_x\phi(x(t|x_0,m(\cdot)))\cdot\oloklirosi{U}{}{f(x(t|x_0,u))}{m_t(u)} =\frac{d}{dt}\phi(x(t|x_0,m(\cdot)))$, for almost all $t\geq 0$. Integration by parts gives $\inner{N^u}{A\phi}{}=\phi(x_0),\,\,\,\,\,\,\mbox{for all}\,\,\,\phi\in\cd$. Equivalently, $A^*N^u=\delta_{x_0}$ and thus $\mathcal{R}_X\subset\mathcal{F}_X$ and $\mathcal{J}(x_0)\le V_X^R(x_0)$.  Next, note that $\mathcal{F}_X=\me^+\cap (A^*)^{-1}(\{\delta_{x_0}\})\subset \{\mu\in\me:\,\,\,\norma{\mu}{}\le\frac{1}{\lambda}\}$. By the Banach-Alaoglu Theorem and weak* continuity of $A^*$, we get that $\mathcal{F}_X\ne\emptyset$ is weakly* compact in $\me$. The solvability of~(\ref{primalp}) follows by the weak* continuity of $\me\ni\mu\rightarrow\inner{\mu}{g}{}$.
\end{proof}
 The dual linear program of~(\ref{primalp}) is
\begin{equation}\label{D.1}
  \mathcal{J}^*(x_0):=\sup_{\phi\in\cd}\{\phi(x_0):\,\,\,A\phi\le g \,\,\,\,\mbox{on}\,\,\,X\times U\}.
\end{equation}
The following notion will be valuable for the interpretation of the feasible solutions of~(\ref{D.1}).
\begin{orismos}
A function $\phi\in\cd$ is called a \emph{smooth subsolution to the HJB equation} if $A\phi\le g$ on $X\times U$.
\end{orismos}
Notice that the dual program~(\ref{D.1}) is always feasible. Indeed, the function $\phi(x):=C$, $x\in X$ with $C\le\frac{-\sup_{(x,u)\in X\times U}|g(x,u)|}{\lambda}$ is feasible for the dual program. Moreover, any smooth subsolution to the HJB equation, gives a global lower bound for the value function in~(\ref{B.7}).
\begin{limma}[Lower bound]\label{globallb}
If $\phi\in\cd$ is a feasible solution for the dual program~(\ref{D.1}), then $\phi\le V_X$ on $X$.
\end{limma}
\begin{proof}Let $x\in X$ s.t. $\mathcal{U}_X(x)\neq \emptyset$. Since $\inner{M^u}{A\phi}{}=\phi(x)$, for all $u(\cdot)\in\mathcal{U}_X(x)$ and $A\phi\le g$, we get $\phi(x)\le \inner{M^u}{g}{}$, for all $u(\cdot)\in\mathcal{U}_X(x)$. Thus, $\phi(x)\le V_X(x)$.
\end{proof}

\begin{thm}[Equivalence]\label{equivalence}
Under Assumptions~\ref{ass:control:model} and~\ref{norelaxationgap}, $V_X(x_0)=\mathcal{J}(x_0)=\mathcal{J}^*(x_0)$.
\end{thm}
 Note that the assumptions in Theorem~\ref{equivalence} are weaker than those in~\cite[Theorem 4.4]{Gaitsgory}, as indicated in the introduction.

 The proof of Theorem~\ref{equivalence}, is divided in several parts. The first step is to prove that strong duality holds.

\begin{limma}[Absence of duality gap]\label{D.2}
	Under Assumption~\ref{ass:control:model}~\ref{ass:H1} and~\ref{ass:H2}, if~(\ref{primalp}) is feasible, then there is no duality gap, i.e., $\mathcal{J}^*(x_0)=\mathcal{J}(x_0)$.
\end{limma}
\begin{proof} By virtue of~\cite[Th. 3.10]{Anderson}, it suffices to prove that the set $D=\{(A^*\mu,\inner{\mu}{g}{}):\,\,\,\mu\in\me^+\}$ is $\sigma(\cds\times\ar,\cd\times\ar)$-closed~\cite[Def. 8.2.1.]{Narici}. Similarly to~\cite[Theorem 2.3 (ii)]{LasserreHenrion}, this follows by the Banach-Alaoglu Theorem and the weak* continuity of $A^*$.
\end{proof}

Next, we consider the case of unconstrained state space. Let $C_0(\ar^n\times U)$ be the Banach space of bounded continuous functions that vanish at infinity with the sup-norm. Then $\mathcal{M}(\ar^n\times U)$ is isometrically isomorphic to $C_0(\ar^n\times U)^*$ and $\mathcal{M}(\ar^n\times U)\subset C(\ar^n\times U)^*$. Let
\begin{equation}\label{unconnstrained ocp}
V(x_0):=\inf_{u(\cdot)\in\mathcal{U}^0} \oloklirosi{0}{\infty}{e^{-\lambda t}g(x(t|x_0,u(\cdot)),u(t))}{t}.
\end{equation}
For each $u(\cdot)$ in $\mathcal{U}^0$, we define the corresponding discounted occupation measure $M^u\in\mathcal{M}(\ar^n\times U)^+$ by~(\ref{oom}), for all $E\in\mathcal{B}(\ar^n\times U)$. Then, (\ref{poom}) holds for every $l\in C(\ar^n\times U)$. We set $\mathcal{M}:=\{M^u: u(\cdot)\in\mathcal{U}^0(x_0)\}$. Then the OCP~(\ref{unconnstrained ocp}) can be written as $V(x_0)=\inf\{\inner{\mu}{g}{}: \mu\in\mathcal{M}\}$.
 As before, we define the linear operator $B:C^1(\ar^n)\rightarrow C(\ar^n\times U)$ by
$
  (B\phi)(x,u):=\lambda\phi(x)-\nabla_x\phi(x)\cdot f(x,u),
$
for all $(x,u)\in \ar^n\times U,\,\phi\in C^1(\ar^n)$. $B$ is linear, bounded and weakly continuous. Moreover, its adjoint $B^*:C(\ar^n\times U)^*\rightarrow C^1(\ar^n)^*$ is well defined and weakly* continuous. We consider the convex optimization problem
\begin{equation*}
P(x_0):=\inf\{\inner{\mu}{g}{}:\mu\in\mathcal{F}\},
\end{equation*}
where,
$
{\mathcal{F}}:=\{\mu\in\mathcal{M}(\ar^n\times U)^+: B^*\mu=\delta_{x_0}\}.
$
Similarly as before, $\mathcal{M}\subset\mathcal{F}$ and $\mathcal{F}$ is convex and weakly* compact in $\mathcal{M}(\ar^n\times U)$. In particular, the following is true.
\begin{limma}\label{convex closure}
Under Assumptions~\ref{ass:H1} and~\ref{ass:H3}, $\mathcal{F}$ is the weak* convex closure of $\mathcal{M}$.
\end{limma}
\begin{proof}
\ref{ass:H1} and~\ref{ass:H3} imply that the value function $V$ is Lipschitz continuous~\cite[Prop. 2.1]{Bardi}. By Rademacher's Theorem, $V$ is almost everywhere differentiable. Therefore, $B V\le g$ for almost all $x\in\ar^n$ and all $u\in U$. W.l.o.g., we assume that $g\in C_0(\ar^n\times U)$. This assumption can be lifted after the proof of Theorem. By using standard mollification techniques as in~\cite[Lemma 3.2]{Fleming2}, we can construct a sequence of approximate smooth subsolutions to the HJB equation for the unconstrained OCP~(\ref{unconnstrained ocp}). That is, there exist $\{V_n\}_n\subset C^1(\ar^n)$ and $\{\varepsilon_n\}_n\subset [0,\infty)$, such that $\lim_{n\rightarrow\infty}\norma{V-V_n}{\infty}=0$, $B V_n\le g+\varepsilon_n$ on $\ar^n\times U$ and $\lim_{n\rightarrow\infty}\varepsilon_n=0$. We will use this result to prove that $V(x_0)=P(x_0)$. Clearly,
$P(x_0)\le V(x_0)$. Assume for the sake of contradiction that $P(x_0)<V(x_0)$. Then, there exists $\mu\in\mathcal{M}(\ar^n\times U)^+$, s.t. $B^*\mu=\delta_{x_0}$ and $\inner{\mu}{g}{}< V(x_0)$. So, for all $n\in\mathbb{N}$,
$V_n(x_0)=\inner{\delta_{x_0}}{V_n}{}=\inner{B^*\mu}{V_n}{}=\inner{\mu}{B V_n}{}\le\inner{\mu}{g+\varepsilon_n}{}$. By taking $n\rightarrow\infty$, we get $V(x_0)\le\inner{\mu}{g}{}$, which is a contradiction. Therefore, we have proved that $V(x_0)=P(x_0)$. Finally, assume for the sake of contradiction that $\overline{\text{conv}(\mathcal{M})}^*\ne\mathcal{F}$. Let $\nu\in\mathcal{F}/\overline{\text{conv}(\mathcal{M})}^*$. By the separation theorem, there exists $g\in C_0(\ar^n\times U)$  and $s\in\ar$ such that $\inner{M^u}{g}{}\geq s$, for all $u\in\mathcal{U}^0$ and $\inner{\nu}{g}{}<s$. Since a function in $C_0(\ar^n\times U)$ can be uniformly approximated by a Lipschitz continuous function in $C_0(\ar^n\times U)$ \cite[pg. 124]{Fleming2}, we can replace the function $g\in C_0(\ar^n\times U)$ with a new Lipschitz continuous function in $C_0(\ar^n\times U)$ and adjust $s$ so that the separation property still holds. Assume that this specific $g$ is the cost function of the OCP~(\ref{unconnstrained ocp}). Then, $P(x_0)<s\le V(x_0)$, which is a contradiction.
\end{proof}
Next, we consider for each $m(\cdot)\in\mathcal{V}^0$, the corresponding relaxed occupation measure $N^m\in\mathcal{M}(\ar^n\times U)^+$ defined by~(\ref{rom}), for each $E\in\mathcal{B}(\ar^n\times U)$ and satisfying~(\ref{prom}), for each $l\in C(\ar^n\times U)$. Set $\mathcal{R}:=\{N^m: m(\cdot)\in\mathcal{V}^0\}$. It is a known result that $\mathcal{R}$ is the weak* closure of $\mathcal{M}$.
\begin{porisma}\label{representation}
For each $\mu\in\mathcal{F}$, there exists $\Lambda\in\mathcal{P}(\mathcal{R})$ such that
 $\inner{\mu}{l}{}=\oloklirosi{\mathcal{R}}{}{\inner{\gamma}{l}{}}{\Lambda(\gamma)}$, for all $l\in C(\ar^n\times U)$.
\end{porisma}
\begin{proof}
Similar arguments as in~\cite[Cor. 1.4]{Vinter} prove the result for all $l\in C_0(\ar^n\times U)$. Since a bounded continuous function can be approximated pointwise by an increasing sequence of continuous functions with compact support, the application of the monotone convergence theorem completes the proof.
\end{proof}
We will now return to the state constrained OCP~(\ref{B.7}). The following Lemma indicates the relation between all the previously defined sets.
\begin{limma}\label{transition}
We make the convention that each measure $\mu\in\mathcal{M}(X\times U)^+$ is a measure on the whole space $\ar^n\times U$ with $\text{spt}\,\mu\subset X\times U$. The following hold,
\begin{eqnarray*}
% \nonumber to remove numbering (before each equation)
  \mathcal{M}_X &=&\{\mu\in\mathcal{M}: \text{spt}\,\mu\subset X\times U\},  \\
   \mathcal{R}_X &=&\{\mu\in\mathcal{R}: \text{spt}\,\mu\subset X\times U\},  \\
   \mathcal{F}_X &=&\{\mu\in\mathcal{F}: \text{spt}\,\mu\subset X\times U\}.
\end{eqnarray*}
\end{limma}
\begin{proof}
We will prove the inclusion $\supset$ in the first assertion. Let $M^u\in\mathcal{M}$, for some $u(\cdot)\in\mathcal{U}^0$, such that $\text{spt}\,M^u\subset X\times U$. Note that $\norma{M^u}{}=\frac{1}{\lambda}$. Therefore, $\oloklirosi{0}{+\infty}{e^{-\lambda t}}{t}=\frac{1}{\lambda}=M^u(X\times U)=\oloklirosi{0}{+\infty}{e^{-\lambda t}\delta_{(x(t),u(t))}(X\times U)}{t}$. Thus, $(x(t),u(t))\in X\times U$, for almost all $t\geq 0$. So, $u(\cdot)\in\mathcal{U}_X(x_0)$.
\end{proof}
We are now ready to prove Theorem~\ref{equivalence}.
\begin{proof}
By Proposition~\ref{minimizer}, there exists a minimizer $\mu_0\in\mathcal{F}_X$, such that $\inner{\mu_0}{g}{}=\mathcal{J}_X(x_0)\le V_X^R(x_0)$. Then, by Corollary~\ref{representation}, there exists a Borel probability measure $\Lambda$ on $\mathcal{R}$ such that $\inner{\mu_0}{l}{}=\oloklirosi{\mathcal{R}}{}{\inner{\gamma}{l}{}}{\Lambda(\gamma)}$, for all $l\in C(\ar^n\times U)$. We set $l(x,u)=\text{d}(x,X)$, where $d(x,X)$ is the euclidean distance of $x$ from $X$. Since, $\textup{spt}\,\mu_0\subset X\times U$, we have $\inner{\mu_0}{\text{d}(\cdot,X)}{}=0$. Therefore, $\inner{\gamma}{\text{d}(\cdot,X)}{}=0$, for $\Lambda$-almost all $\gamma\in\mathcal{R}$. Equivalently, $\Lambda(\mathcal{R}_X)=1$. So, $\mathcal{J}_X(x_0)=\inner{\mu_0}{g}{}=\oloklirosi{\mathcal{R}_X}{}{\inner{\gamma}{g}{}}{\Lambda(\gamma)}$. We have that $\inner{\mu_0}{g}{}\le\inner{\gamma}{g}{}$, for all $\gamma\in\mathcal{R}_X$. Assume for the sake of contradiction that $\inner{\mu_0}{g}{}<\inner{\gamma}{g}{}$, for all $\gamma\in\mathcal{R}_X$. Then, $\inner{\mu_0}{g}{}<\oloklirosi{\mathcal{R}_X}{}{\inner{\gamma}{g}{}}{\Lambda(\gamma)}$ which is a contradiction. Therefore, there exists $\gamma_0\in\mathcal{R}_X$ such that $\mathcal{J}(x_0)=\inner{\mu_0}{g}{}=\inner{\gamma_0}{g}{}\geq V_X^R(x_0)$. So, we have proved that $\mathcal{J}(x_0)=V_X^R(x_0)$ and that the relaxed admissible control associated to $\gamma_0$ is optimal for the relaxed OCP~(\ref{relaxed OCP}). Finally, under Assumption~\ref{norelaxationgap} we have $V_X(x_0)=V_X^R(x_0)=\mathcal{J}(x_0)$.
\end{proof}
%====================================================================================================================================================================================

\section{Primal-Dual Moment-s.o.s. LMIs, Optimal Control Synthesis and Numerical Example}
Throughout this section we make the following assumptions on the data of the OCP~(\ref{B.7}).
\begin{ass}[Polynomial Data] \label{polynomialdata}   \
\begin{enumerate}[label=\text{(A\arabic*)}, itemsep = -1mm, topsep = -3mm]
  \item \label{AA1} $g\in\ar[x,u]$ and $\{f_k\}_{k=1}^n\subset\ar[x,u]$;
  \item \label{AA3} $X\times U$ is a compact basic semi-algebraic set of the form $X\times U=\{(x,u)\in\ar^n\times\ar^m:\,\,\,q_i(x,u)\geq 0,\,\,\,i=1,\ldots,N\},$ for some $N\in\mathbb{N}$ and $\{q_i\}_{i=1}^{N}\subset\ar[x,u]$;
  \item \label{AA4} There exists $\xi\in\ar[x,u]$ such that $\xi\in\left\{\sigma_0+\athroisma{i}{1}{N}{\sigma_i\,q_i}:\,\,\{\sigma_i\}_{i=0}^{N}\subset\Sigma[x,u]\right\}$ and the level set $\left\{(x,u)\in\ar^n\times\ar^m:\,\,\xi(x,u)\ge 0\right\}$ is compact.
\end{enumerate}
\end{ass}
Assumption~\ref{AA4} is not restrictive. Indeed, if necessary, we can add in the representation of $X\times U$ the polynomial inequality $q_{N+1}(x,u):=K-\norma{(u,x)}{2}^2\geq 0$ , where $K$ is an upper bound of $X\times U$. Then, the new equivalent representation satisfies Assumption~\ref{AA4}.

 We now formulate Lasserre's hierarchy of primal-dual moment-s.o.s. semidefinite relaxations \cite{Lasserre} for our specific problem. Under Assumption~\ref{AA4}, we are able to apply \emph{Putinar's Postivstellensatz}~\cite{Putinar} and characterize in a computationally tractable way the moment-sequences of a finite Borel measure on the compact basi semi-algebraic set $X\times U$. So, we conclude that the primal LP~(\ref{primalp}) is equivalent to an infinite SDP problem. By optimizing over finite truncated sequences of moments, we obtain a hierarchy of finite SDP relaxations. In particular, let $\deg(g)=2v_0$ or $2v_0-1$ and $\deg(q_i)=2v_i$ or $2v_i-1$, for all $i=1,\ldots,N$. For each $r\ge\max_{i=0,\ldots,N}v_i$, we consider the LMI-relaxation of order $r$,
\begin{equation} \label{primalrelaxation}
\mathcal{J}_r(x_0):= \left\{
\begin{array}{cl}
\inf\limits_{z=(z_{\gamma})_{|\gamma|\le 2r}}&L_z(g) \\
\text{s.t.}&L_z(h^{(\alpha)})=b^{(\alpha)},\,\mbox{for all}\,\alpha\in\mathcal{A}_r \\
  & M_r(z)\succeq 0,\\
  &M_{r-v_i}(q_i\,z)\succeq 0,\,\,i=1,\ldots,N,
\end{array} \right.
\end{equation}
where
$h^{(\alpha)}(x,u):=\lambda x^{\alpha}-\nabla x^{\alpha}\cdot f(x,u)$ and $b^{(\alpha)}:={x_0}^\alpha$, for each $\alpha\in\mathbb{N}^n$ and $\mathcal{A}_r:=\{\alpha\in\mathbb{N}^n:\,\,\deg(h^{(\alpha)})\le 2r\}$. The dual of~(\ref{primalrelaxation}) is
\begin{equation} \label{dualrelaxation3}
	\mathcal{J}^*_r(x_0):= \left\{
	\begin{array}{cl}
		\sup\limits_{\lambda_{\alpha,\sigma_i}} & \phi(x_0) \\
		\text{s.t.}& \phi(x)=\sum_{\deg h^{(\alpha)}\le 2r}\lambda_{\alpha}x^{\alpha}, \\
		& g- A\phi\,=\,\sigma_0+\athroisma{i}{1}{n}{\sigma_i\,q_i}, \\
		& \lambda_{\alpha}\in\ar,\,\,\sigma_i\in\Sigma[x,u],\\
		&\deg\sigma_i\,q_i\le 2r,i=0,\ldots,N.
	\end{array} \right.
\end{equation}
We highlight that the dual SDP~(\ref{dualrelaxation3}) is a tightening of the dual linear program~(\ref{D.1}), since its feasible solutions are polynomial subsolutions to the HJB equation. Moreover, note that the s.o.s. representation that $g-A\phi$ is required to satisfy can be expressed via SDP feasibility tests, since the degree of the involved s.o.s. polynomials is fixed.

Note that the number of moments, i.e., the the size of the vector $z=(z_{\gamma})_{|\gamma|\le 2r}$, in the LMI relaxation~(14) of order $r$ is $R:=\binom{n+m+2r}{2r}$. Therefore, for a fixed state space with dimension $n$ and control space with dimension $m$, $R$ grows as $O(r^{n+m})$. If the relaxation order $r$ is fixed, then $R$ grows as $O((n+m)^r)$, that is polynomially in the size of the problem. This means that taking into account the current performance of general-purpose SDP solvers, the moment approach is appealing for small to medium size OCPs. One can exploit structure, (e.g., sparsity) of the specific instance OCP under consideration to improve scalability and overcome this computational limit~\cite[Sec. 4.6]{Lasserre},~\cite{Papachristodoulou}. Moreover, in our work, we have chosen the monomial basis to represent polynomials. However, other bases (e.g., Chebyshev) may be more efficient from a computational point of view~\cite{Korda}. One other possible strategy is to develop alternative positivity certificates, which should be less computationally demanding~\cite{MajumdarII}.
\begin{thm}[{{\cite[Thm. 4.3]{Lasserre},\cite[Lemma 5]{kordamci}}}]\label{LMIconvergence}
 Under Assumption~\ref{polynomialdata} and if~(\ref{primalp}) is feasible, then $\mathcal{J}_r^*(x_0)=\mathcal{J}_r(x_0)\uparrow\mathcal{J}(x_0)=\mathcal{J}^*(x_0)$, as $r\rightarrow\infty$. If in addition, Assumptions~\ref{ass:control:model} and~\ref{norelaxationgap} hold, then $\mathcal{J}_r^*(x_0)=\mathcal{J}_r(x_0)\uparrow V_X(x_0)$, as $r\rightarrow\infty$.
\end{thm}
Some preliminary results regarding the convergence rate, are available in the recent work~\cite{rate}.

We next discuss how the  resulting dual s.o.s. SDP relaxations can be used for the computation of an approximate feedback controller. Numerical examples in~\cite{controlsynthesis} for the finite horizon case, have shown that the optimal solution $\phi^*$ of~(\ref{dualrelaxation3}) approximates well the value function $V(\cdot)$ along optimal trajectories starting from $x_0$, but it gives an unsatisfactory approximation for the other points. This observation was formalized in~\cite[Theorem 1]{conicoptimization}. We would like to enlarge this region and have a good approximation of the value function $V$ on a given set $X_0\subset X$. Assume that $U_X(x_0)\ne\emptyset$, for all $x_0\in X_0$. Let $\mu_0$ be the uniform probability measure on $X_0$ and consider the average value $\overline{V}(\mu_0):=\inner{V}{\mu_0}{}=\oloklirosi{X_0}{}{V(x)}{\mu_0(x)}$. Consider the primal averaged linear program
\begin{equation}\label{averaged2}
 \mathcal{J}(\mu_0):=\inf_{\mu\in\me^+}\{\inner{\mu}{g}{}:\,\,\,\,A^*\mu=\mu_0\},
\end{equation}
with corresponding dual
\begin{equation}\label{averaged3}
 \mathcal{J}^*(\mu_0):=\sup_{\phi\in\cd}\{\inner{\mu_0}{\phi}{}:\,\,\,A\phi\le g \,\,\,\,\mbox{on}\,\,\,X\times U\}.
\end{equation}
Intuitively, the primal averaged LP~(\ref{averaged2}) describes a superposition of a possibly uncountable set of OCPs. Under the assumptions of Theorem~\ref{equivalence} and by linearity, we have $\overline{V}(\mu_0)=\mathcal{J}(\mu_0)=\mathcal{J}^*(\mu_0)$.

Assume that a relaxation of order $r$ of the averaged dual program~(\ref{averaged3}) has an optimal solution $\phi$. By routine modifications in~\cite[Theorem 2]{conicoptimization}, we get that $\phi$ is a good approximation of the value function $V$ along all optimal trajectories which start from $X_0$. A natural candidate for feedback control is derived by the following minimization problem: For fixed $x\in X$
\begin{equation}\label{feedback}
  u^*(x)=\arg\min_{u\in U} (g(x,u)-A\phi(x,u)).
\end{equation}
For each $x\in X$, consider a closed neighborhood $\mathcal{S}_x$ of $x$. Similarly to~\cite{controlsynthesis}, we propose the following heuristic iterative algorithm.
\begin{enumerate} [itemsep = -0mm, topsep = -2mm]
  \item Set $\overline{x}=x(t).$
  \item Solve the semidefinite relaxation of the averaged dual program~(\ref{averaged3}) (where $X_0=\mathcal{S}_{\overline{x}}$).
  \item Apply the control law derived by~(\ref{feedback}) until $x(t)\notin\mathcal{S}_{\overline{x}}$.
  \item Go to step 1.
\end{enumerate}
We conclude with an illustrative numerical example. Consider the nonlinear double integrator \vspace{-1mm}
\begin{equation*}
\begin{array}{lll}
 \inf\limits_{u(\cdot)} &  \int_0^{\infty}  &  \hspace{-5mm} e^{-0.1 t}(x_1(t)^2+x_2(t)^2)\drv t  \\
  \mbox{s.t.} & \dot{x}_1(t) & \hspace{-3mm}=x_2(t)+0.1x_1(t)^3, \\
  &\dot{x}_2(t) &\hspace{-3mm}=-0.3u(t), \\
   &x(0) &\hspace{-3mm}=(0,0.7)^{\textup{T}}, \\
   &x(t) &\hspace{-3mm}\in X:=\{x: \,||x||_2\le 1\},\,\,u(t)\in U:=[-1,1]. \vspace{-1mm}
\end{array}
\end{equation*}
A problem with similar data was considered in~\cite[Section 7.1.]{Korda}. The following table shows the approximation of the optimal value for different relaxation orders of the moment SDP~(\ref{primalrelaxation}). It seems that the sequence $(J_r(x_0))_r$ converges to $2.2043$, as $r\rightarrow\infty$. Note that 2.2043 is a lower bound for $V_X(x_0)$. \vspace{-5mm}

\begin{figure}[t]	
	\centering
	\scalebox{1}{\input{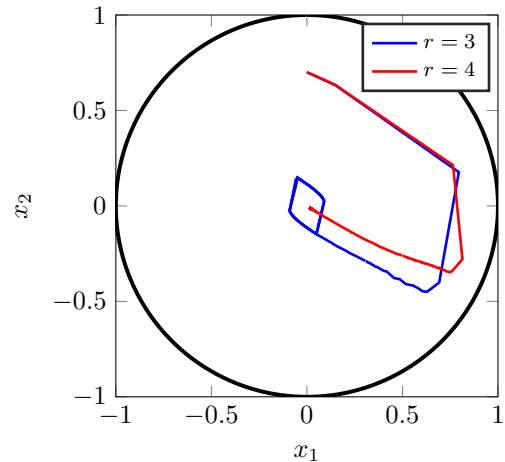}}
	\caption{Trajectories obtained from relaxation order $r=3$ and $r=4$.}
	 \label{fig:traj1}
\end{figure}

\begin{table}[!htb]
\centering
\caption{Numerical results for varying relaxation order $r$. }
\label{tab:Ex1}

  \begin{tabular}{c@{\hskip 2mm} | c@{\hskip 2mm} c@{\hskip 2mm} c@{\hskip 2mm} c@{\hskip 2mm} c@{\hskip 2mm} c@{\hskip 2mm} c}
 $r$   &     $2$  & $3$ & $4$ & $5$ & $6$ & $7$ & $11$\\
 \hline
 $J_r(x_0)$ &  0.1121 &  1.6465 & 2.1978  & 2.2042 & 2.2042 &  2.2043 &  2.2043\\
 \hline
 CPU time&  0.79 & 0.92 & 1.24
  & 2.65 & 6.77 & 17.73  &  677.73
\end{tabular}
\end{table}
%\hfill \break
%\hfill \break
%\begin{center}
%  \begin{tabular}{c l}
%
%  relaxarion order r& $J_r(\mu_0)$\\
%  \hline
%  2&  0.1121\\
%  3& 1.6465\\
%  4& 2.1978\\
%  5& 2.2042\\
%  6& 2.2042\\
%  7& 2.2043\\
%  8& 2.2043\\
% 11& 2.2043
%  \end{tabular}
%  \end{center}
%  \hfill \break
%\hfill \break
\vspace{-2mm}
For the design of an approximate feedback control law, we use the averaged LP formulations~(\ref{averaged2}) and ~(\ref{averaged3}), where the initial condition is described by the uniform probability measure $\mu_0$ on the whole state constraint set $X$.  If an optimal solution $\phi_r^*$ of the averaged dual SDP relaxation of order $r$ is given, then by using~(\ref{feedback}) a feedback control is obtained by
$
u_r(x)=\,\,\mbox{sign}\,\, \frac{\partial}{\partial x_2}\phi^*(x),\,\,\,\,x\in X.
$
Let $V^{u_r}$ be its corresponding cost. Then $V^{u_r}$ is an upper bound for $V_X(x_0)$. Figure~\ref{fig:traj1} displays the trajectories obtained by the dual SDP relaxations of order $r=3$ and $r=4$ . The corresponding costs are $V^{u_3}=2.2479$ and $V^{u_4}=2.2582$. Since the OCP of our example has no analytical solution, we evaluate the performance of the extracted feedback controller by evaluating the gap $G_r=100\frac{V^{u_r}-2.2043}{2.2043}$. We have $G_3=1.98\%$ and $G_4=2.45\%$. This is a strong indication that the extracted feedback controllers are near optimal. In fact, one can verify that in the presented example, Assumptions~\ref{ass:control:model},~\ref{norelaxationgap} and~\ref{polynomialdata} are satisfied (input-affine polynomials dynamics with convex control space). Therefore, Theorems~\ref{equivalence} and~\ref{LMIconvergence} hold. \vspace{-2mm}

\vspace{-2mm}

%\begin{figure}[t]	
%	\centering
%	\scalebox{1}{\input{figures/Ang_figure1.tex}}
%	\caption{Trajectory obtained from relaxation order $r=3$.}
%	 \label{fig:traj1}
%\end{figure}	
%
% Figure~\ref{fig:traj2} displays the trajectory obtained by the dual SDP relaxation of order $r=4$. The corresponding cost is $V^{u_4}=2.2582$.
%\begin{figure}[t]	
%	\centering
%	\scalebox{1}{\input{figures/Ang_figure2_r4.tex}}
%	\caption{Trajectory obtained from relaxation order $r=4$.}
%	 \label{fig:traj2}
%\end{figure}	

%===============================================================================

%===============================================================================

%\begin{figure}[h!]
%\centering
%\includegraphics[width=0.5\textwidth]{figures/ApproxUnifRegularized.eps}
%\caption{Regularized with $\gamma = 1$.}
%\label{fig:UnifApproxDensityRegularised}
%\end{figure}

\bibliographystyle{IEEEtran}
%\bibliography{IEEEabrv,bibliography}
%
%\bibliographystyle{siam}
\bibliography{ref}

% Generated by IEEEtran.bst, version: 1.13 (2008/09/30)
\begin{thebibliography}{10}
\providecommand{\url}[1]{#1}
\csname url@samestyle\endcsname
\providecommand{\newblock}{\relax}
\providecommand{\bibinfo}[2]{#2}
\providecommand{\BIBentrySTDinterwordspacing}{\spaceskip=0pt\relax}
\providecommand{\BIBentryALTinterwordstretchfactor}{4}
\providecommand{\BIBentryALTinterwordspacing}{\spaceskip=\fontdimen2\font plus
\BIBentryALTinterwordstretchfactor\fontdimen3\font minus
  \fontdimen4\font\relax}
\providecommand{\BIBforeignlanguage}[2]{{%
\expandafter\ifx\csname l@#1\endcsname\relax
\typeout{** WARNING: IEEEtran.bst: No hyphenation pattern has been}%
\typeout{** loaded for the language `#1'. Using the pattern for}%
\typeout{** the default language instead.}%
\else
\language=\csname l@#1\endcsname
\fi
#2}}
\providecommand{\BIBdecl}{\relax}
\BIBdecl

\bibitem{Bardi}
M.~Bardi and I.~Capuzzo-Dolcetta, \emph{Optimal control and viscosity solutions
  of {H}amilton-{J}acobi-{B}ellman equations}.\hskip 1em plus 0.5em minus
  0.4em\relax Birkh\"auser Boston, Inc., Boston, MA, 1997.

\bibitem{shooting}
H.~J. Pesch, ``A practical guide to the solution of real-life optimal control
  problems,'' \emph{Control and Cybernetics}, vol.~23, no. 1/2, pp. 7--60,
  1994.

\bibitem{MPC}
\BIBentryALTinterwordspacing
L.~Gr\"une and J.~Pannek, \emph{Nonlinear model predictive control: Theory and
  Algorithms}, ser. Communications and Control Engineering Series.\hskip 1em
  plus 0.5em minus 0.4em\relax Cham, Switzerland: Springer, 2017, second
  edition. [Online]. Available:
  \url{http://dx.doi.org/10.1007/978-3-319-46024-6}
\BIBentrySTDinterwordspacing

\bibitem{direct}
O.~von Stryk and R.~Bulirsch, ``Direct and indirect methods for trajectory
  optimization,'' \emph{Annals of Operations Research}, vol.~37, no.~1, pp.
  357--373, 1992.

\bibitem{NN}
Y.~Tassa and T.~Erez, ``Least squares solutions of the {HJB} equation with
  neural network value-function approximators,'' \emph{Trans. Neur. Netw.},
  vol.~18, no.~4, pp. 1031--1041, Jul. 2007.

\bibitem{HLL2}
D.~Hern\'andez-Hern\'andez, O.~Hern\'andez-Lerma, and M.~Taksar, ``The linear
  programming approach to deterministic optimal control problems,'' \emph{Appl.
  Math. (Warsaw)}, vol.~24, no.~1, pp. 17--33, 1996.

\bibitem{Gaitsgory}
V.~Gaitsgory and M.~Quincampoix, ``Linear programming approach to deterministic
  infinite horizon optimal control problems with discounting,'' \emph{SIAM J.
  Control Optim.}, vol.~48, no.~4, pp. 2480--2512, 2009.

\bibitem{Vinter}
R.~Vinter, ``Convex duality and nonlinear optimal control,'' \emph{SIAM J.
  Control Optim.}, vol.~31, no.~2, pp. 518--538, 1993.

\bibitem{BhattBorkar}
A.~G. Bhatt and V.~S. Borkar, ``Occupation measures for controlled {M}arkov
  processes: characterization and optimality,'' \emph{Ann. Probab.}, vol.~24,
  no.~3, pp. 1531--1562, 1996.

\bibitem{FlemingVermes}
W.~H. Fleming and D.~Vermes, ``Convex duality approach to the optimal control
  of diffusions,'' \emph{SIAM J. Control Optim.}, vol.~27, no.~5, pp.
  1136--1155, 1989.

\bibitem{ref:Peyman-17}
\BIBentryALTinterwordspacing
P.~{Mohajerin Esfahani}, T.~{Sutter}, D.~{Kuhn}, and J.~{Lygeros}, ``{From
  Infinite to Finite Programs: Explicit Error Bounds with Applications to
  Approximate Dynamic Programming},'' \emph{https://arxiv.org/abs/1701.06379},
  Jan. 2017. [Online]. Available: \url{https://arxiv.org/abs/1701.06379}
\BIBentrySTDinterwordspacing

\bibitem{LasserreHenrion}
\BIBentryALTinterwordspacing
J.~B. Lasserre, D.~Henrion, C.~Prieur, and E.~Tr\'elat, ``Nonlinear optimal
  control via occupation measures and {LMI}-relaxations,'' \emph{SIAM J.
  Control Optim.}, vol.~47, no.~4, pp. 1643--1666, 2008. [Online]. Available:
  \url{http://dx.doi.org/10.1137/070685051}
\BIBentrySTDinterwordspacing

\bibitem{controlsynthesis}
D.~Henrion, J.~B. Lasserre, and C.~Savorgnan, ``Nonlinear optimal control
  synthesis via occupation measures,'' in \emph{2008 47th IEEE Conference on
  Decision and Control}, Dec 2008, pp. 4749--4754.

\bibitem{Korda}
M.~Korda, D.~Henrion, and C.~N. Jones, ``Controller design and value function
  approximation for nonlinear dynamical systems,'' \emph{Automatica J. IFAC},
  vol.~67, pp. 54--66, 2016.

\bibitem{FW}
H.~Frankowska and F.~Rampazzo, ``Filippov's and {F}ilippov-–{W}a\.{z}ewski's
  theorems on closed domains,'' \emph{Journal of Differential Equations}, vol.
  161, no.~2, pp. 449--478, 2000.

\bibitem{Gaitsgory2}
V.~Gaitsgory, ``Averaging and near viability of singularly perturbed control
  systems,'' \emph{Journal of Convex Analysis}, vol.~13, no.~2, pp. 329--352,
  2006.

\bibitem{Anderson}
E.~J. Anderson and P.~Nash, \emph{Linear programming in infinite-dimensional
  spaces}.\hskip 1em plus 0.5em minus 0.4em\relax John Wiley \& Sons, Ltd.,
  Chichester, 1987.

\bibitem{Narici}
L.~Narici and E.~Beckenstein, \emph{Topological vector spaces}, 2nd~ed.\hskip
  1em plus 0.5em minus 0.4em\relax CRC Press, Boca Raton, FL, 2011.

\bibitem{Fleming2}
W.~H. Fleming and D.~Vermes, \emph{Generalized Solutions in the Optimal Control
  of Diffusions}.\hskip 1em plus 0.5em minus 0.4em\relax New York, NY: Springer
  New York, 1988, pp. 119--127.

\bibitem{Lasserre}
J.~B. Lasserre, \emph{Moments, {P}ositive {P}olynomials and {T}heir
  {A}pplications}.\hskip 1em plus 0.5em minus 0.4em\relax London: Imperial
  College Press, 2010.

\bibitem{Putinar}
M.~Putinar, ``Positive polynomials on compact semi-algebraic sets,''
  \emph{Indiana Univ. Math. J.}, vol.~42, no.~3, pp. 969--984, 1993.

\bibitem{Papachristodoulou}
Y.~Zheng, G.~Fantuzzi, and A.~Papachristodoulou, ``Exploiting sparsity in the
  coefficient matching conditions in sum-of-squares programming using admm,''
  \emph{IEEE Control Systems Letters}, vol.~1, no.~1, pp. 80--85, July 2017.

\bibitem{MajumdarII}
A.~A. Ahmadi and A.~Majumdar, ``{DSOS} and {SDSOS} optimization: {LP} and
  {SOCP}-based alternatives to sum of squares optimization,'' in \emph{2014
  48th Annual Conference on Information Sciences and Systems (CISS)}, March
  2014, pp. 1--5.

\bibitem{kordamci}
M.~Korda, D.~Henrion, and C.~N. Jones, ``{Convex computation of the maximum
  controlled invariant set for polynomial control systems},'' \emph{{SIAM
  Journal on Control and Optimization}}, vol.~52, no.~5, pp. pp.2944--2969,
  Oct. 2014.

\bibitem{rate}
------, ``Convergence rates of moment-sum-of-squares hierarchies for optimal
  control problems,'' \emph{Systems {\&} Control Letters}, vol. 100, pp. 1--5,
  2017.

\bibitem{conicoptimization}
D.~Henrion and E.~Pauwels, \emph{Linear Conic Optimization for Nonlinear
  Optimal Control}.\hskip 1em plus 0.5em minus 0.4em\relax Philadelphia, USA:
  SIAM, 2016, pp. 121--134.

\end{thebibliography}

\end{document}